\documentclass[letterpaper, 12pt]{article}
\usepackage{amsmath, amssymb, amsthm}
\usepackage{amsfonts}
\usepackage{epsfig, psfrag, harvard}
\usepackage{graphicx}
\usepackage{setspace}
\usepackage[left=3.5cm,top=3.22cm,right=3.5cm,bottom=3.22cm]{geometry}
\citationmode{abbr}
\def\br{{\bf r}}
\def\ba{{\bf a}}
\def\bX{{\bf X}}

\def\bU{{\bf U}}
\def\bD{{\bf D}}

\def\bV{{\bf V}}
\def\bI{{\bf I}}
\def\by{{\bf y}}
\def\bone{{\bf 1}}

\def\R{{\rm R}}

\def\bbeta{{\bf \beta}} 
\def\btheta{{\bf \theta}}

\def\abs#1{|#1|}

\begin{document}

\title{A note on the group lasso and a sparse group lasso}
\author{
{\sc Jerome Friedman}
\thanks{Dept. of Statistics, Stanford Univ., CA
94305, jhf@stanford.edu}\\
{\sc Trevor Hastie}
\thanks{Depts. of Statistics, and Health, Research \&
  Policy, Stanford Univ., CA
94305, hastie@stanford.edu}\\
and
{\sc Robert Tibshirani}\thanks{Depts. of Health, Research \&
  Policy, and Statistics,
    Stanford Univ, tibs@stanford.edu}
}

\maketitle
\date

\begin{abstract}
We consider the group lasso penalty for the linear model.
We note that the standard algorithm for solving the problem 
assumes that the  model matrices in each group are orthonormal.
 Here we consider a  more general penalty that
blends the lasso ($L_1$) with the group lasso (``two-norm'').
This penalty yields solutions that are sparse at both the group
and individual feature levels. We derive an efficient algorithm for
the resulting convex problem based on coordinate descent. This
algorithm can also be used to solve the general form of the group
lasso, with non-orthonormal model matrices.
\end{abstract}

\medskip

\section{Introduction}

 In this note, we consider the problem of prediction  
using  a linear model.
 Our data consist of  $\bf y$,  a vector of
$N$ observations, and $\bf X$, a $N \times p$ matrix of features.                 

Suppose that the $p$ predictors are divided into $L$ groups, with
$p_\ell$ the number in
group $\ell$. For ease of notation, we use a matrix $\bX_\ell$ to
represent the predictors corresponding to the $\ell$th group, with
 corresponding coefficient vector $\beta_\ell$.
Assume that $\bf y$ and $\bX$ has been centered, that is, all variables have mean zero.

In an elegant paper, \citeasnoun{YL2007} proposed the  group lasso which solves the convex optimization problem
\begin{eqnarray}
\min_{\beta\in \R^p}\left(||\by-\bone-\sum_{\ell=1}^L
\bX_\ell\beta_\ell||_2^2 +\lambda \sum_{\ell=1}^L\sqrt{p_\ell} ||\beta_\ell||_2\right),
\label{LM.grlasdef}
\end{eqnarray}
where the $\sqrt{p_\ell}$ terms accounts for the varying group sizes,
and $||\cdot||_2$ is the Euclidean norm (not squared). 
This procedure acts like the lasso at the group level:
depending on  $\lambda$, an
entire group of predictors may drop out of the model.
In fact if the group sizes are all one, it reduces to the lasso.
\citeasnoun{meier2008} extend the group lasso to logistic regression.

The group lasso does not, however, yield sparsity within a group. That is,
if a group of parameters is non-zero,  they will all be non-zero.
In this note we  propose a more general penalty that yields sparsity at both
the group and individual feature levels, in order to select groups
and predictors within a group.
We also point out  that the  algorithm proposed by \citeasnoun{YL2007} for fitting the group lasso
assumes that the  model matrices in each group are orthonormal.
The algorithm that we provide for our more general criterion
also works for the standard group lasso with non-orthonormal model matrices.

\begin{figure}[hbtp]
\begin{center}
\begin{psfrags}
\psfrag{beta1}{$\beta_1$}
\psfrag{beta2}{$\beta_2$}
  \epsfig{file=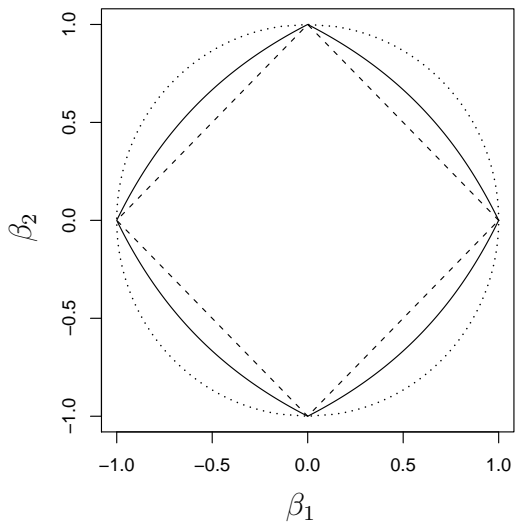,width=.5\textwidth}
\end{psfrags}
\end{center}
  \caption{\small\em Contour lines for the penalty for the group lasso (dotted), lasso (dashed) and sparse group lasso penalty (solid), for a single group
with two predictors.}
  \label{fig:schem}
\end{figure}

We consider the {\em sparse group lasso} criterion
\begin{eqnarray}
\min_{\beta\in \R^p}\left(||\by-\sum_{\ell=1}^L
\bX_\ell\beta_\ell||_2^2 +\lambda_1 \sum_{\ell=1}^L ||\beta_\ell||_2+\lambda_2||\beta||_1\right).
\label{LM.sgrlasdef}
\end{eqnarray}
where $\beta=(\beta_1,\beta_2,\ldots \beta_\ell)$ is the entire  parameter  vector.
For notational simplicity we omit the weights $\sqrt{p_\ell}$.
Expression (\ref{LM.sgrlasdef}) is the sum of convex functions and is therefore convex.
Figure \ref{fig:schem}  shows the constraint region for the
 group lasso, lasso  and sparse group lasso. 
A similar penalty involving both group lasso and lasso terms is discussed
in \citeasnoun{peng2009}.
When $\lambda_2=0$, criterion (\ref{LM.sgrlasdef}) reduces to the group lasso, whose computation we
discuss next.

\section{Computation for the group lasso}
Here we briefly review the  computation for the group lasso of \citeasnoun{YL2007}.
In the process we  clarify a confusing issue regarding orthonormality of
predictors within a group.

The subgradient equations (see e.g. \citeasnoun{Bertsekas1999}) for the group lasso are
\begin{eqnarray}
-\bX_\ell^T(y-\sum_\ell \bX_\ell\beta_\ell)+\lambda\cdot s_\ell =0;\;\ell=1,2,\ldots L,
\end{eqnarray}
where $s_\ell=\beta_\ell/||\beta_\ell||$ if $\beta_\ell \neq 0$ and
$s_\ell$ is a vector with $||s_\ell||_2 <1$ otherwise.
Let the solutions be $\hat\beta_1, \hat\beta_2 \ldots \hat\beta_\ell$.
If
\begin{eqnarray}
||\bX_\ell^T(y-\sum_{k\neq \ell} \bX_k\hat\beta_k)||<\lambda
\label{eqn:zerocond}
\end{eqnarray}
 then   $\hat\beta_\ell$ is zero;
otherwise  it satisfies
\begin{eqnarray}
\hat\beta_\ell=(\bX_\ell^T\bX_\ell +\lambda/||\hat\beta_\ell||)^{-1}\bX_\ell^T r_\ell
\label{two}
\end{eqnarray}
where $$\br_\ell=\by-\sum_{k\neq \ell} \bX_k\hat\beta_k$$

Now if we assume that $\bX_\ell^T \bX_\ell=\bI$, and let
$s_\ell=\bX_\ell^Tr_\ell$,
then (\ref{two}) 
simplifies to $\hat\beta_\ell=(1-\lambda/||s_\ell||) s_\ell$.
This leads to an
algorithm that cycles through the groups $k$, and is a blockwise coordinate
descent procedure. It is given in \citeasnoun{YL2007}.

If however the predictors are not orthonormal,
one approach is to  orthonormalize them before applying the
group lasso.
However this will not generally provide a solution to the original problem.
In detail, if  $\bX_\ell=\bU\bD\bV^T$, then the columns of $\bU=\bX_\ell \bV\bD^{-1}$   are orthonormal.  Then $\bX_\ell \bbeta_\ell = \bU\bV\bD^{-1}
\bbeta_\ell= \bU[\bV\bD^{-1}\bbeta_\ell]=  \bU\bbeta_\ell*$.
But $||\bbeta_\ell*||= ||\bbeta_\ell||$ only if  $\bD= \bI$. This will not be true in general, e.g.
if $\bX$ is a set of dummy varables for a factor, this is true only if the
number of observations in each category is equal.

Hence an alternative approach is needed.
In the non-orthonormal case,
we can think of equation (\ref{two}) as a ridge regression,  with  the ridge parameter
depending on $||\hat\beta_\ell||$.
A complicated scalar equation can be derived  for  $||\hat\beta_\ell||$ from (\ref{two});
then substituting into the right-hand side of (\ref{two})
gives the solution. However this is not a good approach numerically, as
it can involve dividing by the norm of a vector that is very close to zero.
It is also not guaranteed to converge.
In the next section we provide a  better solution to this problem, and to
the sparse group lasso.
\section{Computation for the sparse group lasso}
The criterion (\ref{LM.grlasdef})
is  separable  so that block coordinate descent can be used
for its optimization.
Therefore we  focus on just one group $\ell$, and denote the predictors by 
$\bX_\ell=Z=(Z_1, Z_2, \ldots Z_k)$, the coefficients by $\beta_\ell=\btheta=(\theta_1, \theta_2,
\ldots \theta_k)$ and the residual by $\br=\by-\sum_{k \neq \ell} X_k\beta_k$.
The subgradient equations are
\begin{eqnarray}
-Z_j^T(\br-\sum_j Z_j\theta_j)+\lambda_1 s_j
 +\lambda_2 t_j=0
\label{subgrad}
\end{eqnarray}
for $j=1,2,\ldots k$
 where $s_j=\theta_j/||\btheta||$ if $\btheta_\ell \neq 0$ and $s$ is a
 vector satisfying
$||s||_2 \leq 1$ otherwise, and $t_j \in {\rm sign(\theta_j)}$,
that is $t_j=sign(\theta_j)$ if $\theta_j \neq 0$ and
$t_j \in [-1, 1]$ if $\theta_j = 0$.
Letting $\ba=\bX_\ell\br$, then a necessary  and sufficient condition
for $\btheta$ to be zero is that the 
system of equations
\begin{eqnarray}
a_j=\lambda_1 s_j +\lambda_2 t_j
\end{eqnarray}
have a solution with $||s||_2\leq 1$ and $t_j\in [-1,1]$.
We can determine  this by minimizing
\begin{eqnarray}
J(t)=(1/\lambda_1^2)\sum_{j=1}^k (a_j -\lambda_2 t_j)^2=\sum_{j=1}^k s_j^2
\label{zerocond}
\end{eqnarray}
with respect to the $t_j \in [-1,1]$
and then checking if $J(\hat t)\leq 1$.
The minimizer is easily seen to be
\begin{equation*}
\hat t_j= 
\begin{cases} \frac{a_j}{\lambda_2} & \text{if $\abs{ \frac{a_j}{\lambda_2}}\leq 1$,}
\\
{\rm sign}(\frac{a_j}{\lambda_2}) & \text{if $\abs{ \frac{a_j}{\lambda_2}} > 1$.}
\end{cases}
\end{equation*}

Now if $J(\hat t)>1$, then we must minimize the criterion 
\begin{eqnarray}
\frac{1}{2}\sum_{i=1}^N\Bigl(r_i-\sum_{j=1}^k Z_{ij}\theta_j\Bigr)^2 +\lambda_1 ||\btheta||_2+\lambda_2\sum_{j=1}^k |\theta_j|
\label{three}
\end{eqnarray}
This is the sum of a convex differentiable function (first two terms) and a separable
penalty, and hence we can use coordinate descent to obtain the global
minimum.

Here are the details of the coordinate descent procedure.
For each $j$ let $\br_j=\br-\sum_{k\neq j} Z_k\hat\theta_k$.
Then $\hat\theta_j=0$ if $|Z_j^T\br_j|<\lambda_2$.
This follows easily by examining the subgradient equation corresponding to (\ref{three}).
Otherwise if $|Z_j^T\br_j|\geq\lambda_2$ we  minimize (\ref{three}) by a  one-dimensional
search over $\theta_j$.
We use the {\tt optimize} function in the R package, which
 is a combination of golden section search and
     successive parabolic interpolation.

This leads to  the following algorithm:

\bigskip

\centerline{\underline{\bf Algorithm for the sparse group lasso}}

\begin{enumerate}
\item Start with $\hat\beta=\beta_0$ 
\item In group $\ell$ define  $\br_\ell=\by-\sum_{k \neq \ell} \bX_k\beta_k$,
$\bX_\ell=(Z_1,Z_2 ,\ldots Z_k)$, $\beta_\ell=(\theta_1, \theta_2, \ldots
   \theta_k$) and
 $\br_j=\by'-\sum_{k\neq j} Z_k\theta_k$. Check if $J(\hat t)\leq 1$ according
to (\ref{zerocond}) and if so set $\hat\beta_\ell=0$.
  Otherwise  for $j=1,2,\ldots k$,
if $|Z_j^T\br_j|<\lambda_2$ then $\hat\theta_j=0$;  if instead
$|Z_j^T\br_j|\ge\lambda_2$ then minimize
\begin{eqnarray}
\frac{1}{2}\sum_{i=1}^N(y'_i-\sum_{j=1}^k Z_{ij}\theta_j)^2 +\lambda_1 ||\btheta||_2+\lambda_2\sum_{j=1}^k |\theta_j|
\end{eqnarray}
over $\theta_j$ by a one-dimensional optimization.
\item Iterate step (2) over groups $\ell=1,2,\ldots L$  until convergence.
\end{enumerate}
If $\lambda_2$ is zero, we instead use condition (\ref{eqn:zerocond})
for the group-level test and we don't need to check the condition
$|Z_j^T\br_j|<\lambda_2$.  With these modifications, this algorithm also gives a effective method for solving the
group lasso with non-orthogonal model matrices.

\bigskip

Note that in the special case
where $\bX_\ell^T \bX_\ell=I$, with $\bX_\ell=(Z_1, Z_2, \ldots Z_k)$ then
its is easy to show that
\begin{eqnarray}
\hat\theta_j=\Bigl(||S(Z_j^T\by,\lambda_2)||_2- \lambda_1\Bigr)_+ \frac{S(Z_j^T\by,\lambda_2)}{||S(Z_j^T\by,\lambda_2)||_2}
\label{orth}
\end{eqnarray}
and this reduces to the algorithm of \citeasnoun{YL2007}.

\section{An example}

We generated $n=200$ observations with $p=100$ predictors, in ten blocks of ten.
The second fifty predictors iall have coefficients of zero.
The  number of non-zero coefficients in the first five blocks of 10 are
 (10, 8, 6, 4, 2, 1)  respectively, with coefficients equal to $\pm 1$, the sign
chosen at random.
The predictors are standard Gaussian with correlation 0.2 within a group
and zero otherwise.
Finally, Gaussian noise with standard deviation 4.0 was added to each observation.

\begin{figure}[hbtp] 
  \epsfig{file=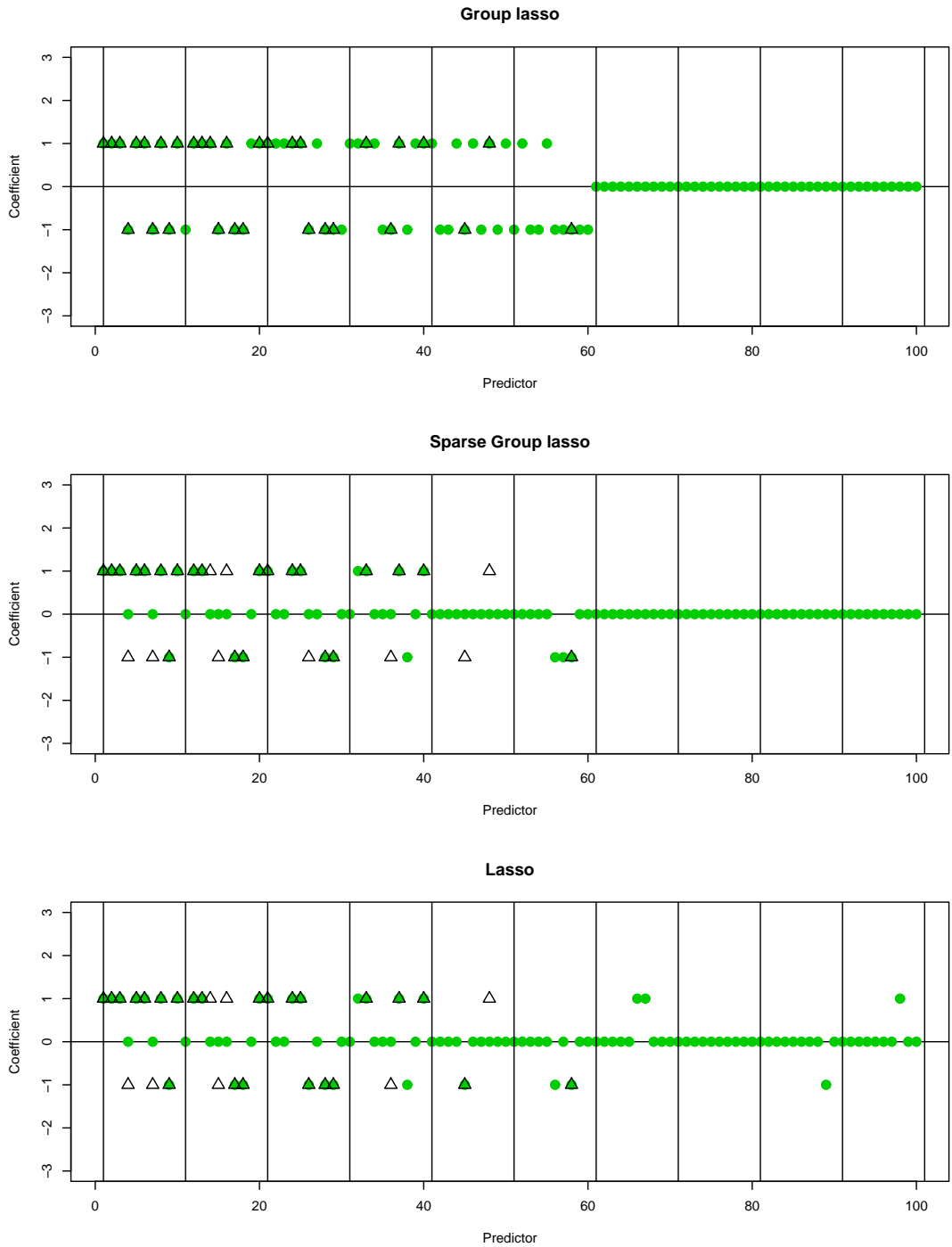,width=\textwidth}
  \caption{\small\em Results for the simulated example. True coefficients
are indicated by the open triangles while the filled green circles
indicate the sign of the estimated coefficients from each method.}
  \label{fig:six}
\end{figure}

\begin{figure}[hbtp]
  \epsfig{file=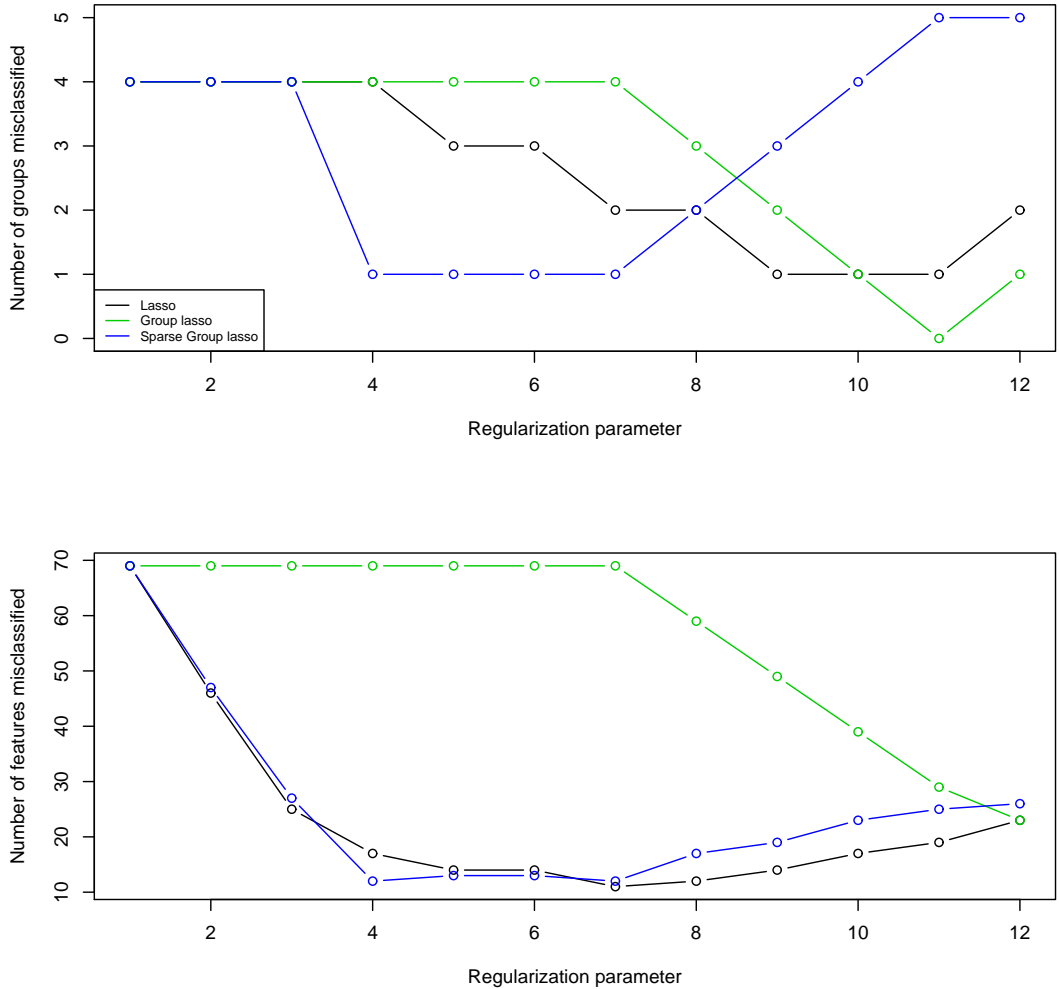,width=\textwidth}
  \caption{\small\em Results for the simulated example.
The top panel shows the number of groups that are misclassified as the
regularization parameter is varied. A misclassified group is one with at
least one nonzero coefficient whose estimated
coefficients are all set to zero, or  vice versa.
The bottom panel shows the number of individual coefficients that
are misclassified, that is, estimated  to be zero when the true coefficient
is nonzero or vice-versa.}
  \label{fig:five}
\end{figure}

Figure \ref{fig:six}  shows the signs of the estimated coefficients from the lasso, group
lasso and sparse group lasso, using a well chosen tuning parameter
for each method (we set  $\lambda_1=\lambda_2$ for the sparse group lasso).
The corresponding misclassification 
rates for the groups and individual features are shown in Figure
\ref{fig:five}. We see that the sparse group lasso strikes an
effective compromise between the lasso and group lasso, yielding
sparseness at the group and individual predictor levels.

\bibliographystyle{agsm}
\bibliography{/home/tibs/texlib/tibs,/home/hastie/bibtex/tibs}

\end{document}